
\magnification=1200
\NoBlackBoxes

\centerline{\bf Singular traces and compact operators. I.}
\vskip1.5truecm
\centerline{S. Albeverio$^{1,4,5}$, D. Guido$^2$,
A. Ponosov$^{1,6}$, S.Scarlatti$^3$}
\vskip2.truecm\noindent
$^1$: Fakult\"at f\"ur Mathematik, Ruhr-Universit\"at Bochum,
Germany\par\noindent
$^2$: Dipartimento di Matematica, Universit\`a di Roma Tor Vergata,
Italy
\par\noindent $^3$: Dipartimento di Matematica, Universit\`a di
L'Aquila, Italy
\par\noindent $^4$: SFB-237, Essen-Bochum-D\"usseldorf
\par\noindent $^5$: CERFIM, Locarno
\par\noindent $^6$: Supported by the Deutsche Forschungsgemeinschaft
\vskip2.truecm\noindent
{\bf Abstract:} We give a necessary and sufficient condition on a
positive compact operator $T$ for the existence of a singular trace
(i.e. a  trace vanishing on the finite rank operators) which takes
a finite non-zero  value on $T$. This generalizes previous results
by Dixmier and Varga. We also  give an explicit description of
these traces and associated ergodic states on $\ell^\infty(\Bbb N)$
using tools of non standard analysis in an essential way.

\beginsection{Section 1. Introduction}\par
\bigskip
In 1966,  Dixmier proved
that there exist  traces on $B(H)$, the bounded linear operators on a
separable complex Hilbert space, which are not normal [D2]. Dixmier traces have
the further property to be ``singular", i.e. they  vanish on
the finite rank operators.
\par
 The importance of this type of traces is well-known due to their
applications in non commutative geometry and quantum field theory (see
[C]).
\par
Since every trace can be decomposed in a sum of a normal and a
singular part,
and since the normal traces on $B(H)$ coincide with the usual trace up
to a constant, it is natural to investigate the structure of the class
of the
singular traces. In this paper we study this problem at a
local level, i.e. we fix an operator $T$ in $B(H)$, and ask the
following
questions:
\item{$(i)$} Do there exist  singular traces which are non-trivial
 ``on $T$", i.e. on the two-sided ideal generated by $T$ in $B(H)$?
\item{$(ii)$} Is it possible to give an explicit
description of these traces, when they exist?
\par
We shall prove that the answer to the first question is positive if
and only if the operator $T$ is generalized eccentric (see definition
2.7).
\par
In relation to the second question, we describe explicitly some
classes of singular traces, which may possibly be extended to larger ideals,
and
we give a detailed analysis of the structure of such classes.
\par
As already mentioned, the first idea  for  constructing singular
traces
is due to  Dixmier, who
considered  compact operators for which the sum of the first $n$
singular values diverges at a given
suitable rate. The singular trace is then obtained evaluating on these
partial sums, appropriately rescaled, a
state on $\ell^\infty (\Bbb N)$ which is invariant under ``2-dilations".
\par
 More recently, the general question which operator ideals in $B(H)$
support
 traces has been  studied by Varga [V].
The procedure used by Varga in order to describe traces differs from
that of
Dixmier in the choice of the states on $\ell^\infty(\Bbb N)$.
 \par
In two preceding papers [AGPS1], [AGPS2] we gave explicit formulas
for Dixmier-type traces and introduced  a new class of singular
traces. The result concerning generalized eccentric operators has been
announced in [AGPS3].
\par
The present paper is organized as follows.
In Section 2 we introduce the basic definition of generalized
eccentric operators and prove the main theorem about singular
traceability. We remark that while the question of mere traceability has a
trivial answer when restricted to the trace class operators,
this is not the case for  singular traceability. Therefore our analysis turns
out to be an  extension of the theory of Varga. Moreover, we illustrate two
different techniques to construct singular traces. These techniques are
a generalization of those in [D2] and [V].
\par
In Section 3 we describe ergodic states giving rise to both kinds of
singular traces introduced in Section 2. The basic technique we use is
related to non standard analysis (NSA).
Section 3 also involves the representation of Banach-Mazur limits by
NSA.
Such representation has been discussed before - e.g. in [KM], [L].
\par
In Section 4 we work out explicitly the computation of the Dixmier traces of an
operator, again using the NSA framework in an  essential way.
\par
Let us finally mention that there are still several problems which deserve
attention. For instance under which conditions a singular trace can be extended
to a larger ideal, and possibly the existence of a maximal ideal in this
context.  Another interesting problem is to find out a general representation
formula for all singular traces. Results on this kind of questions will appear
in
[AGPS4].

\beginsection{Section 2. Singular traces and generalized eccentric
operators}\par
\bigskip
Let $\Cal R$ a von Neumann algebra and $\Cal R^+$ the cone of its
positive
elements. A {\it weight} on $\Cal R$ is a linear map
$$\phi :\Cal R^+\to [0,+\infty ]$$
Any weight can be extended by linearity on the natural domain given by
the
linear span of $\{ T\in \Cal R^+|\phi (T)<+\infty\}$

A weight $\tau$ which has the property:
$$
\tau (T^*T)=\tau (TT^*)\qquad \forall T\in\Cal R
$$
is called a {\it trace} on $\Cal R.$\par
The natural domain of a trace $\tau$ is a two--sided ideal denoted by
$\Cal I_\tau.$ For instance the natural
domains of the trivial traces on $\Cal R$ given by $\tau\equiv0$ and
$\tau\equiv +\infty$ are respectively the ideals $\Cal R$ and $\{ 0\}$
while
the usual trace on $B(H),$
the bounded linear operators on a complex, separable Hilbert space
$H$,
is associated with the ideal $L^1 (H)$ of the trace class
operators.\par
A weight $\phi$ on $\Cal R$ is called {\it normal} if for every
monotonically
increasing generalized sequence $\{ T_\alpha ,\alpha\in I\}$ of
elements of
$\Cal R^+$ such that $T=\sup_\alpha T_\alpha$ one has
$$\phi (T)=\underset \alpha \to{\lim }\phi (T_\alpha )$$
{}From now on the von Neumann algebra $\Cal R$ will be fixed to be
$B(H).$\par
A classical result [D1] concerning normal traces on $B(H)$ is the
following:
\bigskip\noindent
{\bf 2.1 Theorem.} \quad  Every non trivial normal trace on
$B(H)$ is proportional to the usual trace.\par
\bigskip\noindent
By a theorem of Calkin (see [GK]), each proper two-sided ideal in $B(H)$
contains the
finite rank operators and is contained in the ideal $K(H)$ of the
compact
linear operators on $H$. Therefore all traces on $B(H)$ live on the
compact
operators, and the following definition makes sense:
\bigskip\noindent
{\bf 2.2 Definition.}\quad
A trace $\tau$ on $B(H)$ will be call {\sl singular} if it vanishes on
the set
$F(H)$ of finite rank operators.\medskip
\bigskip\noindent
{\bf 2.3 Proposition.}\quad Any trace $\tau$ on $K(H)$ can be uniquely
decomposed
as $\tau = \tau_1+\tau_2,$ where $\tau_1$ is a normal trace and
$\tau_2$ is
a singular trace.
\medskip
\demo {Proof} If $\tau (F(H))\equiv 0$ then the result is obvious. Let
us
suppose there exists $A\in F(H)^+$ such that $\tau(A)=1.$\par
Since $A=\sum\limits_1^N \lambda_i E_i ,$ where $\{ E_i\}$ is a set of
rank
one projectors, $\tau (E_{i_0})= C >0$ for some $i_0.$\par
Since all rank one projectors are unitarily equivalent then $\tau
(E)=C$ for
each rank one projector $E.$ As a consequence $\tau =C\ tr$ $(tr(\cdot
)$
denoting
the usual trace) on rank one projectors and therefore, by linearity,
on  all
$F(H).$\par
Let us set $\tau_2\equiv \tau - C\ tr$ and $\tau _1\equiv C\ tr,$ then
$\tau_2
(F(H))\equiv 0.$ It remains only to show the positivity of $\tau_2.$
For
$A\in K(H)^+$ there exists a sequence $\{ A_n\}$ of finite rank
positive
operators such that $A=l.u.b. A_n.$ For this sequence we have $\tau
(A_n) =
tr (A_n)$ hence $\tau (A)\ge l.u.b.\ \tau (A_n)=C\ (l.u.b.\ tr\
(A_n))=C\ tr
(A),$ since $tr(\cdot )$ is normal. From the previous inequality we
get
$\tau_2 (A)\ge 0.$\hfill $\square$\par

\medskip
In view of this proposition, in the rest of the paper we shall
restrict our
attention to the singular traces.
\medskip
For $T$ a compact operator on $H$, $\{\mu_n (T) \}^\infty_{n = 1}$
will denote the non increasing sequence of the eigenvalues of $| T |$
with multiplicity.\par
We shall also set $\sigma_n (T) \equiv  \sum\limits^n_{k=1} \mu_r
(T).$
\bigskip\noindent
{\bf 2.4 Definition.}\quad Let $T$ be a compact operator. We call
{\it integral sequence} of $T$ the sequence
$\{S_n (T) \}^\infty_{n =0}$ which is an indefinite integral (w.r.t.
the
counting measure) of $\{\mu_n (T) \}^\infty_{n=1},$ i.e. \ $S_n (T) -
S_{n-1} (T) = \mu_n (T), n \geq 1,$ \ and such that
$$S_0 (T) \equiv \cases 0 \qquad \quad T \notin L^1 (H)\\
-tr (T) \quad T \in L^1 (H) \endcases$$
Notice that if $T \notin L^1 (H), \quad S_n (T) = \sigma_n (T), \ n
\geq 1,$
while if $T \in L^1 (H),$ then $S_n (T)=\sigma_n(T)-tr(T)\to 0$ as
$n \to \infty$.
\medskip
\bigskip\noindent
{\bf{2.5} Remark.}\quad If $T$ does not belong to $L^1(H)$ and $\tau$
is a
trace which is finite and non-zero on $| T |$ then $\tau$ is
necessarily singular,
that is, the existence  of traces which are non trivial on $T$ is
equivalent
with the existence of non trivial singular traces on $T.$ Since for
$T \in L^1 (H)$ the existence of a non trivial trace is obvious, it
follows that the relevant question is not the mere ``traceability''
of a compact
operator $T,$ but the existence of a singular  trace which is non
trivial
on  $|T|$.
\par
Let us also notice that a trace $\tau$ is finite on $| T |$ if and
only if
the principal ideal $\Cal I (T),$ i.e. the (two--sided) ideal
generated by
$T$ in $B(H)$, is contained in $\Cal I_\tau.$
\bigskip\noindent
{\bf {2.6} Lemma.}\quad Let $T$ be a compact operator. The following
are equivalent:
\medskip
\item {$(i)$} $1$ is a limit point of the sequence $\Bigl \{\frac
{S_{2 n} (T)}
{S_n (T)} \Bigl \}^\infty_{n=0}$
\medskip
\item {$(ii)$} There exists an increasing sequence of natural numbers
$\{p_k \}$ such that
$\underset{k\to\infty}\to{\lim}\ \frac{S_{kp_{k}}(T)}{S_{p_{k}}(T)}=1$
\medskip
\demo {Proof} $(ii)\Rightarrow(i)$ is obvious.\smallskip
$(i)\Rightarrow(ii)$\par
First we exploit the concavity of the sequence $S_n:=S_n(T)$:
$$
{k-2\over k-1} S_n+{1\over k-1} S_{kn}\leq S_{2n}\qquad
n\in\Bbb N,\quad k\geq 2.
$$
{}From that, with simple manipulations, we get
$$
\left|1-{S_{kn}\over S_n}\right|\leq(k-1)
\left|1-{S_{2n}\over S_n}\right| \qquad k,n\in\Bbb N.
$$
By hypothesis $(i)$, for each $k\in\Bbb N$ there exists $p_k\in\Bbb N$
such that
$$
\left|1-{S_{2{p_k}}\over S_{p_k}}\right|\leq {1\over k^2}
$$
Therefore
$$
\left|1-{S_{k{p_k}}\over S_{p_k}}\right|\leq {k-1\over k^2}
$$
and the thesis follows. \hfill$\square$
\bigskip\noindent
{\bf {2.7} Definition.}\quad A compact operator $T$ which satisfies
one of the equivalent properties of Lemma 2.6 will be called {\it
generalized
eccentric}.
\bigskip\noindent
{\bf {2.8} Remark.}\quad The class of generalized eccentric operators
which
are not in $L^1 (H)$ coincides with the class of eccentric operator
considered
in [V].
\medskip
Let $T$ be a compact operator. Then it is clear that $\Cal I \ (T) =
\bigcup\limits^{+ \infty}_{r=1} \Cal I_r \ (T)$
where $\Cal I_r \ (T)$ is the set of all bounded operators of the form
$$A = \sum^r_{i=1} \ X_i \ T \ Y_i, \qquad X_1, \cdots X_r, Y_1,
\cdots
Y_r \in B (H).$$
The estimates below will be crucial for the rest of the section.
\bigskip\noindent
{\bf {2.9} Proposition.}
\item{$(i)$} An operator $A$ belongs to $\Cal I_r(T)$ if and only if
$$
\exists K\in\Bbb R : \mu_{r(n-1) +1}(A)\leq K\mu_n(T),\qquad n\in\Bbb
N
$$
 \item{$(ii)$} Given $A$, $B$ positive compact operators, then
$$
\sigma_n (A+B)\leq\sigma_n (A)+\sigma_n(B)\leq\sigma_{2n}(A+B),
\qquad n\in\Bbb N $$
\medskip
\demo {Proof} See [V] for $(i)$. See [GK] or [S] for
$(ii)$.\hfill$\square$
\medskip
We can now state and prove the main result of this section.
\bigskip\noindent
{\bf {2.10} Theorem.}\quad Let $T$ be a compact operator. Then the
following are
equivalent:
\item {$(a)$} There exists a singular trace $\tau$ such that $0 < \tau
(| T | ) < + \infty$.
\medskip
\item {$(b)$} $T$ is generalized eccentric.
\medskip
Let us remark that the condition $(i)$ in Lemma 2.6   gives some information on
the rate of convergence of the sequence $\{\mu_n(T)\}$. For instance, notice
that if $\mu_n(T)\sim n^\alpha$ as $\to\infty$, condition $(i)$ implies
 $\alpha =-1$. As a consequence, the only Macaev ideal (see e.g. [GK]) in the
domain of a singular trace is $L^{1^+}(H).$ Other natural examples of
sequences $\{\mu_n\}$ satisfying condition $(i)$  are
$\mu_n\equiv\frac{(\log n)^\alpha}{n}$. Hence, by Theorem 2.10, it follows that
if $T\in K(H)$ and $\mu_n(T)= \frac{(\log n)^\alpha}{n}$ there exists a
singular trace $\tau$  such that $\tau (T)=1.$
When $\alpha <-1,$ the domain of the associated singular traces is contained
in $L^1(H)$ (cf. Remark 2.14).
\medskip
\demo {Proof} $(a)\Rightarrow(b)$
\par
As we observed in Remark 2.5, when $T \notin L^1 (H)$
the existence of a non trivial singular trace on $| T |$
is equivalent to the existence of a non trivial trace ``tout court''
on $| T |$. Therefore, in this case, the result is given by
theorem 1 in [V].\par
We are left with $T\in L^1(H)$. Let us notice that, in this case, we
have  only
to show that
$$
\sup_n\frac{S_{kn}(T)}{S_n(T)}>\frac{1}{3}\qquad \forall k \in \Bbb N.
$$
Indeed, if $1$ is  not a limit point of
$\Bigl\{\frac{S_{2n}(T)}{S_n(T)}\Bigr\}$, then
$\sup_n\frac{S_{2n}(T)}{S_n(T)}=l<1$. As a consequence,
$\frac{S_{2^m n}(T)}{S_n(T)}\leq l^m<{1\over3}$ when $n\in\Bbb N$ and
$m>\Bigl|\frac{\log3}{\log l}\Bigr|$ \medskip
The proof will be given by contradiction, that is we assume
$T \in K (H) $ and $\sup_n\frac {S_{k n} (T)} {S_n (T)} \leq
\frac {1} {3}$ for some fixed $k\in\Bbb N$ and then we prove that any
singular
trace is trivial on $\Cal I (T)$.
\medskip
Since a trace is trivial on a principal ideal $\Cal I$ {\it iff} it is
trivial
on a positive operator $S$ which generates $\Cal I$, we may consider
the following operator $S$, realized averaging over the eigenvalues of
$T$,
 $$
\mu_n (S)\equiv\cases
\frac{S_{k^{L(n)}}(T)-S_{k^{L(n) -1}}(T)}{k^{L(n)}-k^{L(n) -1}}&n>1\\
      \mu_1 (T) & n=1 \endcases
 $$
where $L(n)$ is the integer defined by $k^{L (n) - 1} < n \leq
k^{L (n)}.$
\medskip
 From $(i)$ of Proposition 2.9 follows easily that $S \in \Cal I_k
(T)$ and
viceversa $T \in \Cal I_k (S)$, therefore  $\Cal I (S) = \Cal I (T)$.
 \medskip
Next we notice that the following eigenvalue estimate holds:
$$\mu_n (S) \geq 2 k \mu_{k (n -1) + 1} (S). \eqno (2.1)$$
Indeed, this follows from the definition of $\mu_n (S)$ and the
inequality
$$S_{k^{L (n)}}(T)-S_{k^{L(n)-1}}(T) \geq 2 (S_{k^{L (n) +1}}
(T) - S_{k^{ L (n)}} (T))$$
which is a consequence of the assumptions on $T$.
\medskip
By means of a $k$-dilation procedure, we now construct another
compact positive operator $\tilde S$ such that $\tau (S) = \tau
(\tilde
S)$ for any trace $\tau.$
\medskip
We fix an orthonormal basis of $H$ and describe the operators which
are diagonal w.r.t. this basis by means of the corresponding
eigenvalue sequences
$$\aligned
S &\equiv \{ \mu_1, \mu_2 \cdots \}\\
S_1 & \equiv \{ \mu_1,\underset {(k -1)\ \text {times}}\to{\underbrace
{0 \cdots 0}},
\mu_2,
\underset{(k -1)\ \text {times}}\to{\underbrace {0 \cdots 0 }},
\dots\\
S_2 & \equiv 0, \mu_1 \underset{(k -2)\ \text {times}}\to{\underbrace
{0 \cdots 0 }},
0, \mu_2,
\underset{(k -2)\ \text {times}}\to{\underbrace {0 \cdots 0 }},
\dots\\
&................................................................\\
S_k &\equiv
\underset{(k -1) \text {times}}\to{\underbrace {0 \cdots 0 }}, \mu_1,
\underset{(k -1)\ \text {times}}\to{\underbrace {0 \cdots 0 }}, \mu_2,
\dots\\
\endaligned$$
Then we define
$$
{\tilde S} = \frac {1} {k} \sum^k_{i=1} S_i
$$
By linearity and unitary invariance $\tau (\tilde S) = \tau (S)$
for each trace $\tau.$  Moreover, by construction,
$$\mu_{k (n-1) + j} (\tilde S) = \frac {1}{k} \ \mu_n (S) \quad
\ \forall n \in \Bbb N, \quad j = 1, \cdots, k.$$
Hence,
$$\mu_{k (n -1) + j} (\tilde S) \ge 2 \ \mu_{k (n-1) +1} (S) \ge
2 \ \mu_{k (n - 1) + j} (S) \eqno (2.2)$$
by (2.1).
\medskip
It is evident that (2.2) implies
$$\tau (\tilde S) \geq 2 \tau (S) = 2 \tau (\tilde S)$$
which is impossible if $\tau (S)$ if finite non zero.
\medskip
The proof of $(b)\Rightarrow(a)$ follows immediately by Theorem 2.11
(see below).
\hfill $\square$
\medskip
Let us now discuss possible procedures to construct singular traces on
$K(H).$\par
Our first step is a  generalization of a method suggested in [V], in
order to
built up singular traces $\tau$ associated with generalized eccentric
operators
$T$. To this aim it is useful to introduce a triple $\Omega =
(T, \varphi, \{n_k\})$, where $T$ is generalized eccentric, $\varphi$
is a state
on  $\ell^\infty (\Bbb N)$ which vanishes on $c_0$, the space of
infinitesimal sequences, and $n_k= k \ p_k, k \in \Bbb N$, where
$\{p_k \}$
is the sequence of natural numbers given in Lemma 2.6.
\medskip
With such a triple $\Omega$ we associate a functional $\tau_\Omega$
on the positive part of the ideal $\Cal I (T):$
$$\tau_\Omega (A) \equiv \varphi \left ( \left \{
\frac {S_{n_{k}} (A)} {S_{n_{k}} (T)} \right \} \right )
\quad, \quad A \in \Cal I (T)^+ \eqno (2.3)$$
\bigskip\noindent
{\bf 2.11 Theorem.}\quad Let $T$ be a generalized eccentric operator.
The
functional $\tau_\Omega $ defined in (2.3) extends linearly to a
singular trace
on the ideal $\Cal I (T)$ \medskip
\demo {Proof} The positivity, homogeneity and unitary invariance of
$\tau_\Omega$
are obvious. It suffices to check additivity on positive elements.
\medskip
We take $C \in \Cal I (T)^+.$ Then $C \in \Cal I_r (T)$ for some
$r\geq2$
 and, by Proposition $2.9(i)$, if $2\leq r\leq k$
$$\aligned S_{2 k p_{k}} (C) - S_{k p_{k}} (C) &= \sigma_{2 k p_{k}}
(C) - \sigma_{k p_{k}} (C) \leq\\
& \leq \sigma_{r k p_{k}} (C) - \sigma_{r p_{k}} (C)=\\
&=\sum^{rkp_k}_{j = rp_k+1} \mu_j (C) =
\sum^{k p_k}_{j=p_k+1}\sum^r_{i=1} \mu_{r (j - 1) + i}(C)\leq\\
&\leq\sum^{k p_k}_{j=p_k+1}Kr\mu_j(T)=Kr(S_{k p_{k}} (T) -
S_{p_{k}} (T)) \endaligned$$
As a consequence,
$$\left | \frac {S_{2 n_{k}} (C) - S_{n_k} (C)} {S_{n_{k}} (T)}
\right | \ \leq \ K \ r \ \left | 1 - \frac {S_{p_{k}}
(T)} {S_{k p_{k}} (T)} \right |
\underset {k \to \infty} \to {\longrightarrow} 0.$$
Now suppose $T\in L^1(H)$ and $A,\ B\in\Cal I(T)^+$.
Then $A, B \in L^1 (H),$ and Proposition $2.9 (ii)$ implies
$$\frac {S_{n_{k}} (A+B)} {S_{n_{k}} (T)} \geq
\frac {S_{n_{k}} (A)} {S_{n_{k}} (T)}+
\frac {S_{n_{k}} (B)} {S_{n_{k}} (T)}\ge
\frac {S_{ 2 n_{k}} (A+B)} {S_{n_{k}} (T)}$$
Since $\varphi$ is positive and vanishes on infinitesimal sequences,
$$\aligned \tau_\Omega (A+B) & \geq \tau_\Omega (A) + \tau_\Omega (B)
\geq \tau_\Omega (A+B) + \varphi
\left ( \left \{ \frac {S_{2n_{k}} (A+B) - S_{n_{k}} (A+B)}
{S_{n_{k}} (T)} \right \} \right )  =\\
& = \tau_\Omega (A+B) \endaligned$$
i.e. $\tau_\Omega $ is additive.\par
The proof for $T \notin L^1 (H) $ is analogous. \hfill $\square$
\medskip The singular traces $\tau_\Omega$ given by (2.3) which are
associated
with generalized eccentric operators $T$ give rise to a constructive
proof of
the  implication $(b)\Rightarrow(a)$ of theorem 2.10.\par
We shall call such traces {\it generalized Varga traces}. Indeed, in
the case $T \notin L^1 (H),$ the traces $\tau_\Omega$ correspond to
traces constructed in [V].\par
Let us notice that the traces given by (2.3) can be written also as
$$\tau_\Omega (A) = \varphi^{\{n_k\}} \Bigl ( \Bigl \{ \frac {S_{n}
(A)}
{S_{n} (T)} \Bigr \} \Bigr ), \quad A \in \Cal I (T) \eqno (2.4)$$
where $\varphi^{\{n_k\}}$ is the (non normal) state on
$\ell^\infty (\Bbb N)$ defined by
$$\varphi^{\{n_k \}} ( \{a_n \}) \equiv \varphi ( \{ a_{n_{k}} \})$$
We remark that if $\varphi$ is an extremal state on $\ell^\infty (\Bbb
N),$
so is $\varphi^{\{ n_k \}}.$\medskip
Let us now consider in general the functional
$$\tau_{\psi} (A)\equiv \psi \left ( \left \{ \frac {S_{n} (A)}
{S_{n} (T)} \right \} \right )\quad , \quad A \in \Cal I (T) \eqno
(2.5)$$
where $\psi$ is a generic state on $\ell^\infty (\Bbb N).$\par
The above remarks show that, if $\psi$ is chosen as the state
$\varphi^{\{n_k \}},$ then (2.5) gives rise to a singular trace.
\medskip
Other singular traces can be obtained by choosing states $\psi$ in
(2.5) with suitable invariance properties. Generalizing an idea of
Dixmier [D2], we shall prove the following theorem.
\bigskip\noindent
{\bf 2.12 Theorem.}\quad If $\psi$ is a two-dilation invariant state
and
$\underset {n\in\Bbb N} \to {\lim} \ \frac {S_{2_{n}} (T)}
{S_{n} (T)} = 1,$ then
$\tau_{\psi} $ is a trace on $\Cal I (T).$
Moreover, in
this case formula (2.5) gives rise to a singular trace (which will be
denoted by
$\tau_{\psi}$ as well) even on the (larger) ideal
 $${\Cal I}_m(T) \equiv \left
 \{ A \in K(H) {}\vert \left \{ \ \frac {S_n (A)} {S_n (T)}
\right \} \in \ell^\infty \right \}. $$
\medskip
We would like to point out that, when $T\not\in L^1(H)$, the ideal
 ${\Cal I}_m(T)$ is a maximal norm ideal in the sense of Schatten
[S] (see also [GK]).
\bigskip\noindent
{\bf 2.13 Remark.}\quad A 2-dilation invariant state is necessarily
not normal,
more precisely it is zero on the space $c_0$ of infinitesimal
sequences.
Indeed if $\{a_n\}$ has only a finite number of non zero elements then
$\varphi (\{a_n\})=\varphi (\{a_{2^kn}\})=\varphi (\{0\})=0$ for a
sufficiently large $k.$\par
By continuity this result extends to $c_0.$\par
\medskip
\demo {Proof of Theorem 2.12} Unitary invariance is obvious. We prove
positivity. If $ A, T\in L^1(H)$ or $A, T\notin L^1(H)$ then
$\frac{S_n(A)}{S_n(T)}\ge 0$ and positivity follows.\par
For $A\in L^1(H), T\notin L^1(H),$ we have $\frac{S_n(A)}{S_n(T)}\to
0$
as $n\to\infty $, so $\tau_\psi\left(\left\{ \frac{S_n(A)}{S_n(T)}
\right\}\right)=0 $ by Remark 2.13.\par
Finally if $A\notin L^1(H), T\in L^1(H)$ then $A\notin \Cal I (T)$ and
therefore $\tau_\psi (A)=+\infty .$ Now we prove linearity.\par
First we observe that if $A, B$ are such that $A+B\notin \Cal I (T)$
then at least one of them, say $A,$ does not belong to $\Cal I (T).$
Therefore $\tau_\psi (A)+\tau_\psi (A)+\tau_\psi (B).$\par
{}From Proposition $2.9(ii)$ it follows
$$S_n (A+B)\le S_n (A)+S_n (B)\le S_{2n}(A+B)$$
when $A, B\in L^1 (H)$ or $A, B\notin L^1(H).$ In both cases if
$T\notin L^1 (H)$ we get
 $$
\psi\left(\left\{\frac{S_n(A+B)}{S_n (T)}\right\}\right)
\le\psi\left(\left\{\frac{S_n(A)}{S_n (T)}\right\}\right)+
\psi\left(\left\{\frac{S_n(B)}{S_n (T)}\right\}\right)\le
\psi\left(\left\{\frac{S_{2n}(A+B)}{S_n (T)}\right\}\right)
\eqno (2.6)
$$
while if $T\in L^1 (H)$ we get the reversed inequalities.\par
Let us now remark the following property of the state $\psi :$
if $\{b_n\}\in l^\infty (\Bbb N)$ and $\frac{a_n}{b_n}\to 1$ as
$n\to\infty $ then $\psi \left(\{ a_n\}\right)=\psi\left(
\{ b_n\}\right).$ This follows from the fact that $\{a_n -b_n\}
\in c_0.$ Applying these properties we get
$$\aligned
\psi\left(\left\{\frac{S_{2n}(A+B)}{S_n(T)}\right\}\right) &=
\psi\left(\left\{\frac{S_{2n}(A+B)}{S_n(T)}\right\}\left\{
\frac{S_{2n}(T)}{S_n(T)}\right\}\right)=\psi\left(\left
\{\frac{S_{2n}(A+B)}{S_{2n}(T)}\right\}\right)\\
&= \psi\left(\left\{\frac{S_{n}(A+B)}{S_n(T)}\right\}\right)
\endaligned\eqno (2.7)$$
the last equality being a consequence of 2--dilation invariance.\par
Therefore by (2.6), (2.7) and the definition of $\tau_\psi$ we
get
$$\tau_\psi (A+B)=\tau_\psi (A)+\tau_\psi (B)$$
when $A, B\in L^1(H)$ (or $A, B\notin L^1 (H))$ and either
$T\in L^1(H)$ or not.\par
The case $A\in L^1(H), B\notin L^1(H)$ and $T\in L^1(H)$ implies
$A+B\notin L^1(H)$ and therefore $A+B\notin \Cal I (T),$ a situation
already discussed. It remains the possibility $A\in L^1(H), B
\notin L^1(H)$ and $T\notin L^1(H).$ In such a case we have
$$\frac{S_{n}(A+B)}{S_n(T)}\le \frac{S_n(A)+tr (A)+S_n(B)}{S_n(T)}
\le \frac{S_{2n}(A+B)}{S_{2n}(T)}$$
and therefore, since $\frac{tr (A)}{S_n(T)}\to 0$ as $n\to\infty ,$
we obtain that linearity holds once again. This ends the proof.
\hfill$\square$\par
\bigskip\noindent
{\bf 2.14 Remark.}\quad
We notice that if the operator
$T\notin L^1(H)$, then the traces described in theorem 2.12 are exactly
the traces discussed by Dixmier in [D2]. Indeed the sequence
$\{S_n(T)\}$ has all the
properties of the sequence $\{\alpha_n\}$ listed in the paper of
Dixmier. On the other hand given any sequence $\{\alpha_n\}$
with the properties required by Dixmier there exists a generalized
eccentric
operator $T$ for which $S_n (T)=\alpha_n$ (see e.g. [AGPS 1]).\par
In the case $T\in L^1(H)$ our theorem 2.12 produces a new class of
non normal traces, which is in a sense the inverse image inside $L^1(H)$
of the class of Dixmier traces. For such a reason we shall call {\it
generalized Dixmier traces} the singular traces given by Theorem 2.12.
\par
The existence of this new type of traces was announced in [AGPS 2].
\bigskip\noindent
{\bf 2.15 Remark.}\quad
According to the decomposition in Proposition 2.3, a trace is non
normal {\it iff} $\tau_2$ is non-zero, and it is faithful {\it iff}
$\tau_1$ is non-zero. On the other hand if $\tau_2$ vanishes on
$L^1(H)$ it gives no contribution to the sum. Therefore the traces
which come from generalized eccentric operators inside $L^1(H)$,
summed with the usual trace, give the first example of non-normal,
faithful traces.

\def\Re{\Bbb R}
\def\Na{\Bbb N}
\def\R*{{^*\Re}}
\def\N*{{^*\Na}}
\def\claim(#1){\bigskip\noindent{\bf #1}\quad}

\beginsection{Section 3. Two-dilation invariant states and
ergodicity}\par
\bigskip
The main problem we are going to discuss in this section concerns
extremal (ergodic) states which give rise to singular traces. In
our opinion it is non standard analysis (NSA) which supplies the
most convenient tools for this purpose.
\par
Recall that if $\{a_n\}_{n\in\Bbb N}$ a standard sequence of real
numbers,
$\{{}^*a_n\}_{n\in{}^*\Bbb N}$ will denote its non standard
extension.\par
As always, if $x\in{}^*\Bbb R$ then ${}^\circ x\in \Bbb R$ will denote
the
standard part of $x.$\par

We first briefly discuss
extremal states corresponding to the generalized Varga traces
described in Theorem 2.11.
\medskip
Let us denote by $\Delta_{\{n_k\}}$ the set of all extremal points
in the set of non normal states of the form (2.4).
\claim( 3.1 Proposition. ) The set $\Delta_{\{n_k\}}$ consists of
the states
$$\psi(\{a_n\}) \equiv {}^{{}^\circ} ({}^* a_{n_m})$$
for some $m \in {}^* \Bbb N_\infty = {}^* \Bbb N - \Bbb N$.
\demo {Proof} It immediately follows from the fact that extremal
states on $\ell^\infty$ are all of the form
$$\varphi(\{a_n\}) = {}^{{}^\circ} ({}^* a_{m}), {} m \in {}^* \Bbb
N. \eqno(3.1)$$ Since additionally $\varphi$ must vanish on the set
$c_0$ then
$m$ becomes infinitely large. \hfill $\square$ \medskip
\claim( 3.2 Remark.) Of course, instead of infinitely large numbers
one can equivalently use the Stone-{\v C}ech compactification
$\overline\Na$ of $\Na$ and the isomorphism
$\ell^\infty(\Na)\simeq C(\overline \Na)$ given by the
Gelfand transform in order to describe extremal states in
Proposition 3.1. Namely, they will be given by Dirac measures
supported by the set $\overline\Na - \Na$. On the contrary, the
classification of ergodic 2-dilation invariant states {\sl does}
require NSA (see e.g. the  remark 3.6). \medskip
Now we come to the much more difficult problem of classification of
two dilation invariant states. First, we remark that
in order to prove the existence of such states,
Dixmier invoked the amenability of the affine group.
As promised,  we shall adopt here an
alternative point of view, which relies on the use of
NSA and related methods (see e.g. [HL], [AFHKL]).\medskip
\claim(3.3 Theorem.) \quad The map $\omega \to \varphi_\omega$,
$\omega
\in {}^* \Bbb N_\infty$, defined by
$$\varphi_\omega (a) \equiv {}^{{}^\circ} \biggl ( \frac {1} {\omega}
\sum^\omega_{k =1}  {}^* a_{2^k} \biggr ) \eqno (3.2)$$
takes values in the convex set of 2--dilation invariant states
over $\ell^\infty.$
\medskip
\demo {Proof} Let $b_n \equiv \frac {1} {n}
\sum^{n}_{k=1} a_n.$ Since $\{ a_n \}$ is bounded $\{b_n \}$
is also bounded  so that
$\varphi_\omega (a) = {}^{{}^\circ} ({}^* b_\omega) $ is well defined
for
all $\omega.$ Obviously, $\varphi_\omega$ is a state. It is
also 2--dilation invariant since:
$$\varphi_\omega ( \{ a_{2 n} \} ) - \varphi_\omega ( \{a_n \} )
= {}^{{}^\circ} \biggl ( \frac {1} {\omega}
( {}^* a_{2^\omega + 1} - {}^*a_2) \biggr ) = 0.$$
\hfill $\square$ \medskip
A consequence of this theorem is that an explicit formula
for the previously introduced traces can easily be given.
\claim(3.4 Corollary.) If $T$ is an operator verifying
$\lim_n\frac {S_{2n}(T)}{S_n(T)}=1$
 and $\omega$ is an infinite hypernatural number then
$$ \tau_\omega (A) \equiv {}^{{}^\circ} \biggl ( \frac {1} {\omega}
\sum^\omega_{k = 1} \frac {{}^*S_{2^k} (A)} {{}^* S_{2^k} (T)}
\biggr) \qquad  A \in  \Cal I_m(T) \eqno (3.3)$$
is one of the singular traces described in Theorem 2.12.
\medskip
The proof of this corollary follows immediately from Theorems 2.12 and
3.3.
\medskip
There is a simple generalization of the formula (3.2) which describes
2-dilation invariant states. If $j \in {}^* \Bbb N$ and $n
\in {}^* \Bbb N_\infty $ \ the map
$$\{ a_k \} \to {}^{{}^\circ} \biggl ( \frac {1} {n} \sum^n_{i=1} {}^*
a_{j
2^i}
\biggr ) \eqno (3.4)$$
is a 2--dilation invariant state over $\ell^\infty $ and
therefore gives rise to a singular trace.
\medskip
Since any hypernatural $j$ can be written in a unique way as a
product of an odd number and a power of $2,j = (2m -1) 2^{k-1},$
we may rewrite the previous states as
$$\varphi_{k,m,n} (a) = {}^{{}^\circ} \biggl ( \frac {1} {n}
\sum^{k+n}_
{i=k+1} {}^* a_{(2m - 1) 2^{i-1}} \biggr ) \quad k, m \in
{}^* \Bbb N, n \in {}^* \Bbb N_\infty \eqno (3.5)$$
In the rest of this section we shall study states of the form (3.5)  in
relation
to the problem of ergodicity.

Let $\Delta:\Na\to \Na$ be the multiplication by 2, $\Delta_*$ the
corresponding
morphism on $\ell^\infty(\Na)$, $\Delta_*(\{a_n\})=\{a_{2n}\}$, we
shall say
that the
state $\varphi$ is $\Delta$-{\sl invariant} if
$\varphi\circ\Delta_*=\varphi$.
We shall give necessary conditions for extremality in the (convex
compact) set
of $\Delta$-invariant states in terms of NSA.

It is known (see e.g. [E, p.113]) that the states on
$\ell^\infty(\Na)$ can be
identified with the finitely
additive probability measures on $\Na$,
therefore we shall denote any such a state by $\mu$, and the notation
$\mu(A)$
with $A\subset\Na$ makes sense.\par
Moreover extremality of a $\Delta$-invariant state $\mu$ can be
expressed in
terms of ergodicity of $\mu$ seen as  a measure, i.e.
$\mu$ is $\Delta$-invariant if
for each
$A\subset \Na$
 such that $\Delta A=A,$ one has $\mu(A)=0$ or 1.
 \claim(3.5 Remark.) \quad Using
the Stone-\v Cech compactification $\overline\Na$ of
$\Na$ and the isomorphism $\ell^\infty(\Na)\simeq C(\overline \Na)$
given by the
Gelfand transform once again we get an identification of
the states on
$\ell^\infty(\Na)$ with
the $\sigma$-additive probability Radon measures on $\overline\Na$.
On the other hand a transformation on $\Na$ extends to a continuous
transformation on $\overline\Na$. We shall denote with $\overline\mu$,
$\overline\Delta$ the measure and the transformation on $\overline\Na$
induced
by $\mu$ and $\Delta$ respectively. It turns out that ergodicity of
$\overline\mu$ is equivalent to ergodicity of the finitely additive
measure
$\mu$. This equivalence can be shown using well known criteria for
ergodicity
(see e.g. [Ma]).
\bigskip
 \claim(3.6 Remark.) \quad Let us consider the
correspondence $\eta:\Na\times\Na\to\Na$ defined by $(m,n)\to(2m-
1)2^{n-1}$,
which is a bijection. It induces an isomorphism:
$\eta_*:\ell^\infty(\Na)\to \ell^\infty(\Na\times\Na)$ given by
$$(\eta_*a)_{m,n}\equiv a_{\eta(m,n)},$$
$\Delta_*   \equiv    \eta_*^{-1}\Delta\eta_*$    becoming    the
translation $T$ in the second variable: $\Delta_*(m,n)=(m,n+1)\equiv
(m,Tn)$.
It might be thought that
the isomorphism $\eta_*$ gives rise to the splitting of the
dynamical system $(\overline\Na , \overline\Delta )$
in a product
of two dynamical systems $(\overline\Na , id)$ and $(\overline\Na ,
\overline T)$,
thus furnishing a standard
approach to the considerations we shall make below.
Unfortunately this is not true since the spaces
$\overline\Na\times\overline\Na$ and $\overline{\Na\times\Na}$ are not
homeomorphic (see e.g. [G]).
Our idea is to exploit the advantages of NSA, in
particular
the
nice functorial property
$\N*\times\N* = ^*(\Na\times\Na)$.
\medskip
Let $M_\Delta$ denote the set of extremal $\Delta$-invariant (i.e.
ergodic)
states on $\ell^\infty(\Na)$.

\claim(3.7 Proposition.) Any $\mu\in M_\Delta$ coincides with one of
the
states
$\varphi_{k,m,n}$ for some $k,m\in\N*$, $n\in\N*_\infty$.
\medskip
\demo {Proof} Evidently, for each state
$\mu$ and each finite dimensional $E\in \ell^\infty (\Na) \equiv
\ell^\infty$ there exist numbers
$s(E)$ and
$\lambda_j(E)$ such that\medskip
(i) $\mu(a)=\sum^{s(E)}_{j=1}\lambda_j(E) a_j\quad(\forall a\in E),$
\quad(ii) $\lambda_j(E)\ge 0$,
\quad(iii) $\sum^{s(E)}_{j=1} \lambda_j(E)=1$.
\medskip\noindent
For every hyperfinite dimensional space $E$ (see e.g. [AFHKL, p.55]),
$\ell^\infty\subset E\subset{^*{\ell^\infty}}$, there is therefore an
internal
set
$\{\lambda_1,\dots,\lambda_s\}$ satisfying (ii) and (iii) and such
that
$$\mu(a)={^*\mu}({^*a})=\sum^s_{j=1}\lambda_j{^*{a_j}}\qquad(\forall
a\in
\ell^\infty)$$
If $\mu$ is 2-dilation invariant, then for all finite $n$
$$\mu(a)=\sum^s_{j=1}\lambda_j{1\over n}\sum^n_{i=1}{^*a_
{j2^{i-1}}}\qquad(\forall
a\in \ell^\infty).$$
As it can be easily checked, for each finite dimensional space
$F\subset \ell^\infty$,
and each $n\in\Na$ we have
$$\alpha _{n,{^*F}}\equiv\sup_{a\in {^*F}, \Vert a\Vert\le
1}\sup_{1\le
j\le s}
\Biggl\{\vert
{^*\mu}(a) - \sum^s_{j=1}\lambda_j{1\over
n}\sum^n_{i=1}a_{j2^{i-1}} \vert + $$ $$ + \vert
{^*\varphi}_{k,m,n}(a) - {1\over n}\sum^n_{i=1}a_{j2^{i-1}} \vert
\Biggr\}\approx 0 \eqno(3.6)$$
where $j = (2m-1) 2^{k-1}$.
By saturation, there are a  hyperfinite dimensional space $F$ ,
$\ell^\infty\subset F\subset{^*{\ell^\infty}}$ and a number
$n\in\N*_\infty$
such that
(3.6) remains valid. This means that $${^*\mu}(a)\approx
\sum^s_{j=1}\lambda_j
{^*\varphi}_{m,k,n}(a)\quad(j=(2m-1) 2^{k-1})\eqno(3.7)$$
for all $a\in\ell^{\infty}$ and
 immediately implies that the set of all $T-$invariant states
coincides
with  the closed convex hull
$\overline{co}M$ of the set
$M\equiv\{\mu_{k,n}\mid k\in\N*,\, n\in\N*_\infty\}$.

Finally, we show that $M$ is closed which, according to [E, p.708],
would
imply
the inclusion $M_\Delta\subset M$. Consider a directed set
$\Gamma\subset\N*^3$ and assume that
$\varphi_{k,m,n}{\buildrel\Gamma\over\rightarrow}\mu$.    Clearly, for
every
finite
dimensional
$E\subset \ell^{\infty}$ there exists a subsequence
$\{(k_i,m_i,n_i)\}\in\Gamma$
such that
$$^*\varphi_{k_i,m_i,n_i}(a)\approx{^*\mu(a)}\quad{\text for\,
some}\quad
i\in\N*_\infty
\quad{\text and\, for\, all}\quad a\in{^*E},\, \|a\|\le 1.$$
If $\mu$ is $\Delta$ invariant then by saturation
the latter relation holds for a certain hyperfinite dimensional
$E\supset
\ell^\infty$ with $k\in\N*$ and $n\in\N*_\infty$. This concludes the
proof.
\hfill$\square$

\claim(3.8 Corollary.) The problems of description of extremal 2-
dilation
invariant and extremal translation invariant states are equivalent.
\demo{Proof}
Recalling the map $\eta _*$ given by
 $(\eta _* a)_{m,n} = a_{\eta (m,n)}$ where
$\eta_*(m,n)=(2m-1)2^{n-1}$ and applying the proposition just
proved we conclude that for each fixed $m$ the states
$\varphi_{k,m,n}\circ
\eta_*$ coincide with the translation invariant states $\mu_{k,n}$ on
$\ell^{\infty}$ defined by
$$\mu_{k,n}(a) \equiv ^{{}^\circ} \biggl ( \frac {1} {n}
\sum^{k+n}_{i=k+1}
{}^*
a_{i}\biggr ),$$
so that any $\mu \in \eta _*^{-1}(M_{\Delta})$ is contained in
one of the sets
$\{ \delta _m  \otimes \nu \vert \nu \in
M_T\cap M \}$ where $\delta _m$ is given by
$\delta _m (a)\equiv {}^{^\circ} (^* a_{m}),{} m \in\N*$,{} $M_T$
stands for the set of extremal
translation invariant states on $\ell^{\infty}$ and $M\equiv
\{\mu_{k,n}
\vert {} k\in\N*, n\in\N*_{\infty}\}$.
On the other hand it is known (see [KM], [L], or also [AGPS1] where
the
result was proved independently) that $M_T\subset M$ which completes
the
proof. \hfill$\square$

\claim(3.9 Remark.) Of course, the representation given in Proposition
3.7 is not
unique due to Corollary 3.8 and the following trivial
\claim(3.10 Proposition.) If ${k-l\over n}\approx 0$ and ${p-n\over
n}\approx
0$,
then $\mu_{k,n} = \mu_{l,p}.$ \medskip
Now we formulate the main result in this section more precisely.
\claim(3.11 Theorem.) If $\mu$ is an extremal 2-dilation invariant
state
on $\ell^\infty$ then $\mu=\varphi_{k,m,n}$
(with $\varphi_{k,m,n}$ defined by (3.5)) for some $m\in\N*$ and
infinitely
large hypernaturals $k$
and $n$ such that ${n\over k}\approx 0$.\medskip
\medskip
\demo{Proof}
By virtue of Corollary 3.8 it suffices to show that if $\nu\in M_T$
then
$\nu=\mu_{k,n}$ for some
$k, n\in \N*_\infty$. We first prove that $k\in\N*_\infty$. Suppose it
is
not
the case and $k$ is finite. Without loss of generality we can assume
$k=1$,
and, due to proposition 3.10, $n=2m$. If we show that
$\mu_{1,n}\not=\mu_{1,m}$
then the representation $\mu_{1,n}={1\over 2}(\mu_{1,m}+\mu_{m,m})$
implies
$\mu_{1,n}$ is not extremal.

For $b\equiv m2^{-p}$, choose  $p\in\N*_\infty$ such that
$0<{^ob}<\infty$ and
consider a non decreasing sequence $\{b_j\}$ such that
$^*b_{2p}\approx b$.
Define a sequence $\{c_j\}$ of natural numbers by putting
$c_j\equiv[2^jb_j]$
($[\cdot]$ denotes the integer part).
Since
$$c_j\le{1\over 2}[2^{j+1}b_j]+1\le{1\over 2}[2^{j+1}b_{j+1}]+1=
{1\over 2}c_{j+1} +1,$$
one gets $c_{j+1}\ge 2c_j-2$. At the same time,
$$\eqalign{^o\left({^*c_{2p}-m\over m}\right)
&={^o\left({^*c_{2p}-m\over 2^{2p}}\right)}
{^o\left({2^{2p}\over m}\right)}={^o\left({[2^{2p}{^*b_{2p}}]\over
2^{2p}}-b\right)}{1\over {}^ob}\cr
&\le{1\over{^ob}}{^o({^*b_{2p}}-b)}=0.\cr}$$
By the same reason, ${^*c_{2p+1}-n\over n}\approx 0$. Applying
proposition 3.10
once again we obtain
$${1\over
m}\sum^m_{i=1}{^*a_i}\approx{1\over{^*c_{2p}}}\sum^{^*c_{2p}}_{i=1}{^*
a_i}\,;
\quad {1\over n}\sum^n_{i=1}{^*a_i}\approx{1\over
c_{2p+1}}\sum^{^*c_{2p+1}}_{i=1}{^*a_i}\qquad(\forall a\in
\ell^\infty)\eqno(3.8)$$ Now we introduce a set
$B=\cup^\infty_{q=1}[c_{2q-1},c_{2q}]\subset \Na$ and a sequence
$\chi=\{\chi_i\}$ where $\chi_i=1$ for $i\in B$ and 0 otherwise. By
(3.8),
$$\eqalign{\mu_{1,n}(B)&
={^o\left({1\over{^*c_{2p+1}}}\sum^{^*c_{2p+1}}_{i=1}{^*\chi_i}\right)
}=\cr
&{^o\left({{^*c_{2p}}\over{^*c_{2p+1}}}\right)}{^o\left({1\over{^*c_{2
p}}}
\sum^{^*c_{2p}}_{i=1}{^*\chi_i}
+{1\over{^*c_{2p}}}\sum^{^*c_{2p+1}}_{i={^*c_{2p}}+1}{^*\chi_i}\right)
}=\cr
&={^o\left({[2^p{^*c_p}]\over[2^{p+1}{^*c_{p+1}}]}\right)}\mu_{1,m}(B)
\le{1\over 2}\mu_{1,n}(B).\cr}$$
It remains to prove that $\mu_{1,m}(B)\not=0$. In order to see this,
let us
observe that $^*\chi_i=1$ for $i\in[{^*c_{2p-1}},{^*c_{2p}}]$ and that
$^*c_{2p}-{^*c_{2p-1}}\ge {^*c_{2p-1}}-2$. Hence,
$\sharp\{i\mid{^*\chi_i}=1\}\ge{1\over 2}{^*c_{2p}-1}$
(where $\sharp$ means cardinality)
and
$$\mu_{1,m}(B)
={^o\left({1\over{^*c_{2p}}}\sum^{^*c_{2p}}_{i=1}{^*\chi_i}\right)}
\ge{^o\left({1\over{^*c_{2p}}}\left({1\over 2}{^*c_{2p}}-
1\right)\right)}
={1\over 2},$$
so that $\mu_{1,n}$ is not extremal.
We continue the proof assuming $k\in\N*_\infty$ and ${n\over
k}\not\approx 0$
or, equivalently, $^o\left({k\over n}\right)<\infty$. We have to show
that
$\mu_{k,n}$ is not extremal. First we notice that
$$0\le{-k+[\sqrt{k}+1]^2\over n}\le{2\sqrt{k}+1\over
n}\approx{2\sqrt{k}\over
n}={2k\over n\sqrt{k}}\approx 0$$
and analogously,
$${-k+[\sqrt{k}]^2\over n}\approx 0,\qquad
{k+n-[\sqrt{k+n}]^2\over n}\approx 0,\qquad
{k+n-[\sqrt{k+n}+1]^2\over n}\approx 0.$$
Applying now proposition 3.10, one can assume that $k=(2r)^2$,
$k+n=(2s)^2$.

Let us introduce a set $C\equiv\cup^\infty_{i=1}[(2i-1)^2,(2i)^2]$ and
observe
that
$$\mu_{k,n}(C\triangle TC)\le{^o\left({\sharp\{i\mid
i^2\in[k,k+n]\}\over
n}\right)}={^o\left({\sharp\{i\mid 2r\le i\le 2s\}\over (2s)^2-
(2r)^2}\right)}
=0.$$
Extremality of $\mu$ would imply, therefore, that $\mu_{k,n}(C)$
should have
been equal to 0 or 1. On the other hand,
$$\mu_{k,n}(C)={^o\left({1\over n}\sum^s_{i=r+1}\left((2i)^2-(2i-
1)^2\right)
\right)}={^o\left({(4s+1+4r+3)(r-s)\over 8r^2-8s^2}\right)}={1\over
2}.$$
This contradiction implies the result.\hfill$\square$
\claim(3.12 Remark.) The necessary conditions in theorem 3.11 are
surely not
sufficient. To see this, one can consider a state $\mu\equiv\mu_{4^p-
n,n}$ for
arbitrary $p,\, n\in \N*_\infty$. If $n\ge 2^{2p-1}$, then non
extremality of
$\mu$ follows from the theorem. Otherwise, we may introduce the set
$A\equiv
\cup^\infty_{q=1}(2^{2q-1},2^{2q}]$ which is easily seen to be $\mu-
a.e.$
$T-$invariant, but $\mu(A)={1\over 2}.$
\medskip
Now we give a corollary which relates the results of this section with
the
description of the generalized Dixmier traces (for the proof cf. [AGPS1]).
 \claim(3.13 Corollary.) Let $\tau$ be a generalized Dixmier trace on the
ideal ${\Cal
I}_m(T)$  (see Theorem 2.12). Then $\tau$ is in the closure of the convex
hull of
the family
 $$\{\tau_{k,m,n}\mid m\in\N*,\,k,n\in\N*_\infty,\, ({n\over
k})\approx 0\},$$
 where $\tau_{k,m,n}=\tau_{\varphi_{k,m,n}}
$ is the trace associated with the state
$\varphi_{k,m,n}$ given by (3.5) via formula (2.5)
on the same  domain ${\Cal I}_m(T)$.
\claim(3.14 Remark.) The states $\mu_{k,n}$ can be looked upon
intuitively as
averages
on
 intervals of the set $\N*$. This suggests to call {\sl ergodic} all the
intervals associated with ergodic states.\par
Then it is easy to show that if the interval $I$ is ergodic and a
subinterval
$J$  is such that ${|I|\over|J|}\not\approx 0$ (where $|I|$ denotes
the length
of $I$), then $\mu_I=\mu_J$. A sketch of the proof is the following:
let
$I=I_0\cup J\cup I_1$ be a partition of $I$ into subintervals. It
turns out
that
$${}^o\left({|I_0|\over|I|}\right)\mu_{I_0}+
{}^o\left({|J|\over|I|}\right)\mu_{J}+
{}^o\left({|I_1|\over|I|}\right)\mu_{I_1}=\mu_I,$$
hence, by the ergodicity of $I$, $\mu_I=\mu_J$.
\beginsection{Section 4. A computational example}\par
\bigskip
We shall now discuss some advantages of representing singular
traces by means of NSA.
\par
A remarkable advantage lies, in our opinion, in the increased
computability of the value of a singular trace on a given operator
when such a trace is parameterized by some infinite number.
\par
In what follows, we shall work out an example in which we
explicitly calculate the value of the Dixmier trace of an
operator, even though it depends on the
non-standard parameter.
\par
To this aim we shall make use of formula (3.3) choosing a compact
operator $T$
such that $S_n(T)=\log n$. The choice of "summing" logarithmic
divergences has extensively been used by Connes in some applications
to
non-commutative geometry [C].
 \par
Let $q \ge 1$ be a fixed natural number,
we consider a positive compact operator $A_q$ whose sequence of
eigenvalues $(\lambda_n \mid n =  3,4. \dots)$ is defined in the
following way: let $(n_k \mid k = 0,1, \dots)$ be an unbounded
increasing
sequence of natural numbers (with $n_0 \equiv 1$) whose explicit
dependence on $q$ will be given below. For  $n \in (2^{n_k},
2^{n_{k+1}}]$, we define
$$\lambda_n := \quad \frac {n_{k+1} - n_k} {2^{n_{k+1}} - 2^{n_k}}
\eqno (4.1)$$
For $m \ge 2$ we consider the sum $\sigma_{2^m} : = \sum \limits
^{2^m}_{j =3} \lambda_j.$
\medskip
Let $n_k < m \le n_{k+1}$, then we have
$$\sigma_{2^m} = n_k + \frac {2^m - 2^{n_k}} {2^{n_{k+1}} -
2^{n_k}} \cdot (n_{k+1} - n_k) - 1 \eqno (4.2)$$
since
$$\sigma_{2^m}= \sum \limits^{k-1}_{r=0} \sum
\limits^{2^{n_{r+1}}}_{j = 2^{n_r}+ 1} \lambda_j +
\sum \limits^{2^m}_{j = 2^{n_k}+ 1} \lambda_j.$$
\par
Now let $p > 1$ and hence $n_s < p \le n_{s+1}$ for some $s,$
we have
$$\frac {1} {p} \sum^{p}_{m=1} \ \frac {\sigma_{2^m}} {\log 2^m} =
\frac {1} {\log 2} \cdot \frac {1} {p} \left (\sum^{s -1}_{k =0} \ \
\sum^
{n_{k+1}}_{m= {n_k + 1}} \frac {\sigma_{2^m}} {m} + \sum^{p}_{m = n_s
+ 1} \frac {\sigma_{2^m}} {m} \right) \eqno (4.3)$$
We now proceed to estimate the sums appearing on the r.h.s. of
(4.3). \medskip
By means of (4.2) we have
$$\aligned \sum^{n_{k+1}}_{m= n_{k} +1} \frac {\sigma_{2^m}}
{m} &= \left [n_k -1 - \frac {2^{n_k}} {2^{n_{k+1}} - 2^{n_k}}
(n_{k+1} - n_k ) \right] \sum^{n_{k+1}}_{m = n_{k} + 1 }
\frac {1} {m}\\
&+ \frac {n_{k+1} - n_k} {2^{n_{k+1}} -2^{n_k}} \cdot \sum^{n_{k+1}}_
{m = n_{k} + 1} \frac {2^m} {m} \endaligned \eqno (4.4)$$
We notice that the following equalities hold
$$\sum^{n_{k+1}}_{m = n_{k} + 1 } \frac {1} {m} =
\log \frac {n_{k+1}} {n_k} + O  \left (\frac {1} {n_k} -
\frac {1} {n_{k+1}} \right) \eqno (4.5a)$$
$$\sum^{n_{k+1}}_{m = n_{k} + 1 }  \frac {2^m} {m}
= \frac {1} {\log2} \left (\frac {2^{n_{k+1}}} {n_{k+1}} -
\frac {2^{n_k}} {n_k} \right) \left ( 1 + O  \left ( \frac
{1} {n_k} \right) \right) \eqno (4.5b), $$
from which it follows
$$\sum^{n_{k+1}}_{m = n_{k} + 1 } \frac {\sigma_{2^m }} {m}=
\left [ n_k \log\frac {n_{k+1}} {n_k}+ O  \left (1-
\frac {n_k} {n_{k+1}} \right) \right] \left[ 1 + O
\left (\frac {1} {n_k}\right)\right] \eqno (4.6)$$
under the assumption $O  \left (\frac 1{n_k} \right)
\ge O  \left (\frac {2^{n_k}}{2^{n_{k+1}}} \right).$
\medskip
To verify such a condition we fix the initial sequence $(n_k |
k = 0,1, \dots )$ to be of the form $n_k := 2^{kq}$, where
$q \in \Bbb N .$\par
Formula (4.6) takes then the form:
$$\sum^{n_{k+1}}_{m = n_{k} + 1 } \frac{\sigma_{2^m}}{m}=
\left [2^{kq} q \log 2  + O  (1) \right] \left[ 1+ O
(2^{-k q}) \right] \eqno (4.7)$$
Therefore we obtain
$$\frac 1{p\log 2}\sum\limits^{s-
1}_{k=0}\sum\limits^{n_{k+1}}_{m={n_k+1}}
\frac {\sigma_{2^m}}{m}=\frac 1{q\log 2}\left(p\log 2\frac{2^{sq}-
1}{2^q-1}
+O(s)\right)\eqno (4.8)$$
Now, by definition, taking $p\in {}^*\Bbb N_\infty ,$ we have
$$\aligned
\tau^{\text{Dix}}_p (A_q) &={}^\circ\left(\frac 1p\sum\limits^p_{m=1}
{}^*\left(\frac{\sigma_{2^m}}{\log 2^m}\right)\right)\\
&={}^\circ\left(\frac 1{p\log 2}\sum\limits^{s-
1}_{k=0}\sum\limits^{n_{k+1}}_
{m=n_k+1}{}^*\left(\frac{\sigma_{2^m}}{m}\right)\right)+\\
&+{}^\circ\left(\frac1{p\log
2}\sum\limits^{p}_{m=n_s+1}{}^*\left(\frac{\sigma_{2^m}}
m\right)\right)\\
&={}^\circ\left(\frac{2^{sq}}p\right)\frac q{2^q-
1}+{}^\circ\left(\frac 1
{p\log
2}\sum\limits^p_{m=n_s+1}{}^*\left(\frac{\sigma_{2^m}}{m}\right)\right
)
\endaligned\eqno (4.9)$$
where, in the last equality, we have used (4.8) and the fact that
${}^\circ \left(\frac s p\right)=0.$\par
To end the computation of the Dixmier trace of $A_q$ we need to
evaluate
the second term on the r.h.s. of (4.9).\par
We have, for $p\in \Bbb N,$
$$\aligned
\frac 1{p\log 2}\sum\limits^p_{m=n_s+1}\frac{\sigma_{2^m}}{m} &=\\
&=\frac 1{\log 2}\left(\frac{n_s}{p}\right)\left(1+O\left(\frac 1{n_s}
\right)\right)\left(\log\frac p{n_s}+O\left(\frac
1{n_s}\right)\right)+\\
&\frac 1{p(\log 2)^2}\left(\frac{n_{s+1}-n_s}{2^{n_{s+1}}-
2^{n_s}}\right)
\left( \frac{2^p}{p}-\frac{2^{n_s}}{n_s}\right)\left(1+O\left(\frac
1{n_s}
\right)\right)
\endaligned\eqno (4.10) $$
Hence, by estimates similar to the previous ones, and taking $p\in{}^*
\Bbb N_\infty $ we obtain
$${}^\circ\left(\frac 1{p\log
2}\sum\limits^{p}_{m=n_s+1}{}^*\left(\frac
{\sigma_{2^m}}m\right)\right)=\frac 1{\log 2}
{}^\circ\left(\frac {2^{sq}}p\right)
\log{}^\circ\left(\frac p{2^{sq}}\right)
\eqno (4.11)$$
{}From (4.9) and (4.11) it follows
$$\tau_p^{Dix}(A_q)=t\left(\frac q{2^q-1}-\log_2(t)\right)\eqno
(4.12)$$
where $t:={}^\circ \left(\frac {2^{sq}}p\right).$\par
In general $t$ can take any value in the interval $\left[2^{-q},
1\right].$
In particular, in the case $p=2^{sq+r}$, $1\le r\le q,$ formula (4.12)
becomes
$$\tau^{Dix}_p (A_q)= 2^{-r}\left( \frac{q}{2^q-1}+r\right)\eqno
(4.13)$$
\bigskip\noindent

\beginsection{References}\par
\bigskip

\item{[AFHKL]} Albeverio, S., Fenstad, J.E., H\o egh-Krohn, R.,
Lindstr\o m,
T.: {\it Non standard methods in stochastic analysis and mathematical
physics.}, Acad. Press, Orlando (1986).

\item{[AGPS1]} { Albeverio S., Guido D., Ponosov A., Scarlatti S.},
``{\it
Dixmier traces and non standard analysis}'', Proceedings of the III
International Conference ``Stochastic Processes:  Physics and
Geometry'',
Locarno June 1991, Edts. S. Albeverio, D. Merlini,
 World Scientific, Singapore (1993).

\item{[AGPS2]} { Albeverio S., Guido D., Ponosov A., Scarlatti S.},
``{\it Non
standard representation of non normal traces}'', Proceedings of the
German-French ZiF-Meeting ``Dynamics of the  Complex and Irregular
Systems''
Bielefeld, Germany, December 1991, Edts. Ph. Blanchard, L. Streit,
D. Testard,
World Scientific, Singapore (1993).

\item{[AGPS3]} { Albeverio S., Guido D., Ponosov A., Scarlatti S.}
``{\it Singular traces and non-standard analysis.}'', to appear in the
Proceedings of the Blaubeuren Conference on Non Standard Analysis - 1992,
Edts. S. Albeverio, W. A. J. Luxemburg, M. Wolff, Kluwer (1993).

\item{[AGPS4]} { Albeverio S., Guido D., Ponosov A., Scarlatti S.}
``{\it Singular traces and compact operators. II.}'', In preparation.

\item{[C]} Connes, A.: {\it G\'eometrie non commutative},
Intereditions, Paris (1990).

\item{[D1]} Dixmier, J.: {\it Von Neumann algebras}, North Holland
(1981)

\item{[D2]} Dixmier, J.: {\it Existence de traces non normales}, C.R.
Acad.
Sci. Paris {\bf 262}, (1966).

\item{[E]} Edwards, R.E.: {\it Functional analysis, theory and
applications},
New York (1965).

\item{[G]} Glicksberg, I.: {\it Stone-\v Cech compactifications of
products},
Trans. AMS {\bf 90} (1959).

\item{[GK]} Gohberg, I., Krein, M.G.: {\it Introduction to the theory
of
non-selfadjoint operators}, Mos- cow (1985).

\item{[HL]} Hurd, A.E, Loeb, P.A.: {\it An introduction to non-
standard real
analysis}, Acad. Press, Orlando (1985).

\item{[KM]} Keller, G., Moore, L.C., Jr.: {\it Invariant means on the
group of integers}, in: Analysis and Geometry, Ed.
B. Fuchssteiner and W.A.J. Luxemburg, Wissenschaftsverlag,
(1992) 1.

\item{[L]} Luxemburg, W.A.J.:
{\it Nonstandard hulls,
generalized limits and almost convergence},
 in: Analysis and Geometry, Ed.
B. Fuchssteiner
and W.A.J. Luxemburg, Wissenschaftsverlag,
(1992) 19.

\item{[Ma]} Man\'e, R.: {\it Ergodic theory and differentiable
dynamics},
Springer, Berlin-Heidelberg-N.Y. (1987).

\item{[S]} Schatten, R.: {\it Norm ideals of completely continuous
operators},
Springer, Berlin-Heidel- berg-N.Y. (1970).

\item{[V]} Varga, J.V.: {\it Traces on irregular ideals}, Proc. of Am.
Math.
Soc. {\bf 107}, (1989) 715.

\end